\documentclass[12pt]{article}

\usepackage{amsfonts,eucal,bm,amsmath,amssymb,amsthm}
\usepackage[cp1251]{inputenc}
\usepackage[russian]{babel}

\newcommand{\Diff}{{\rm{Diff}}}
\newcommand{\diam}{{\rm{diam}}}
\newcommand{\id}{{\rm{id}}}

\newcommand{\varsloj}{{x}}
\newcommand{\SigmaSix}{{\widetilde{\Sigma}^6}}
\newcommand{\SigmaSixS}{{\SigmaSix \times S^1}}

\newenvironment{pf}{ \noindent {\bf Доказательство.} }
    {$\blacksquare$\medskip}

\newcommand{\supp}{\mathop{\rm supp}}

\newtheorem{theorem}{Теорема}
\newtheorem{lemma}{Лемма}

\theoremstyle{definition}
\newtheorem{definition}{Определение}

\author{
    М.Нальский\footnote{
    Московский Государственный Университет им.~М.~В.~Ломоносова,
    Независимый Московский Университет.
    Работа осуществлена при частичной
    поддержке грантов CRDF~RM1-2358-MO-02, РФФИ~02-01-00482,
    РФФИ~02-01-22002 и~050102801-НЦНИЛ-а.}
}

\title{}
\date{}

\begin{document}

\begin{titlepage}

{\center \large

Московский Государственный Университет им.~М.~В.~Ломоносова

Механико-математический факультет

}

\vspace{1.5cm}

{\flushright \large

УДК 517.9

}

\vspace{2.0cm}

{\center \large Нальский Максим Борисович

}

\vspace{1.0cm}

{\center \Large  Негиперболичность инвариантных мер 

на максимальном аттракторе.

}

\vfill

{\center\large 

Москва, 2007

}

\end{titlepage}

\addtocounter{page}{1}
\tableofcontents


\newpage

\section{Основной результат}

Настоящая работа связана с попытками ответить в той или иной мере на вопрос,
на\-сколь\-ко поведение {\it типичной} динамической системы {\it гиперболично}.

Гиперболические динамические системы имеют ненулевые показатели
Ляпунова, но не являются типичными в пространстве всех динамических
систем \cite{AS}. Неравномерно гиперболические
системы, изучаемые в теории Песина \cite{Pe}, также имеют
ненулевые по\-ка\-за\-те\-ли Ляпунова и так же нетипичны.

Теория Песина изучает диффеоморфизмы, имеющие ненулевые показатели
Ляпунова относительно некоторой фиксированной инвариантной меры,
и описывает поведение траекторий, типичных относительно этой меры.
В общем случае, инвариантная мера может быть задана изначально и согласована
с гладкой структурой, а может определяться динамической системой.
Эти два случая существенно отличаются друг от друга.

Для $C^1$-диффеоморфизмов двумерных многообразий, сохраняющих
площадь, Bochi~\cite{BO} обнаружил, что типичные отображения этого
класса являются либо Аносовскими, либо имеют нулевые показатели
Ляпунова.

Для б\'{о}льших размерностей в работе~\cite{BFP} установлено, что
в пространстве устойчиво-эргодических $C^2$-диффеоморфизмов
с $C^1$-топологией, сохраняющих гладкую форму объ\-е\-ма
на компактном многообразии размерности 2 и выше, существует
открытое и плот\-ное подмножество неравномерно гиперболических отображений.

Если диффеоморфизм частично гиперболичен, сохраняет объем и 
является ус\-той\-чи\-во эргодическим, то Baraviera и Bonatti~\cite{BB} показали,
что сколь угодно малым $C^1$-возмущением можно сделать центральный
показатель Ляпунова ненулевым.

Настоящая работа относится ко второму из упомянутых направлений -- 
исследованию систем, чьи инвариантные меры определяются динамикой 
и могут быть не согласованы с гладкой структурой. Оказывается, 
если отказаться от заданной apriori гладкой инвариантной меры,
то нулевые Ляпуновские показатели встречаются гораздо чаще.

\begin{theorem}[Основной результат]\label{main}
Для замкнутого многообразия $M$, $\dim M \geq 4$, най\-дет\-ся такая область 
$U \subset \Diff^1(M)$, что любой диффеоморфизм $f \in U$ имеет 
ло\-каль\-но максимальный частично гиперболический аттрактор $\Lambda \subset M$
и не\-ато\-мар\-ную эр\-го\-ди\-чес\-кую инвариантную меру $\mu$ с $\supp \mu = \Lambda$, 
один из показателей Ляпунова от\-но\-си\-тель\-но которой равен нулю.
\end{theorem}

Доказательство теоремы~\ref{main} состоит из двух частей. 
В первой части проводится редукция общего вопроса для гладких 
динамических систем к аналогичному результату про косые произведения,
следуя работе~\cite{G05}.

В второй части, более технической и оригинальной, содержится обоснование 
упо\-мя\-ну\-то\-го результата для косых произведений над соленоидом
Смейла-Вильямса. В упро\-щен\-ном случае (подковы Смейла
вместо соленоида) эти рассуждения приведены в ра\-бо\-тах~\cite{GIKN} 
и~\cite{KN06}.

Автор признателен А.~С.~Городецкому, Ю.~С.~Ильяшенко и
В.~А.~Клепцыну за полезные идеи и обсуждения.

\section{Редукция гладкого случая к гельдеровым косым произведениям}

Мы рассматриваем динамические системы, являющиеся косыми произведениями 
над соленоидом Смейла-Вильямса:
\begin{equation}
\label{solen_s1}
F:\Lambda_0 \times S^1\to \Lambda_0 \times S^1, \quad F(s,x)=(T(s),f_s(x)).
\end{equation} 
Здесь $\Lambda_0$~--- соленоид, то есть максимальный аттрактор 
соответствующего отображения $T$ полнотория $D$ в себя 
(см. раздел~\ref{SECTION_skew}); $f_s$ непрерывно зависящее
от параметра $s\in\Lambda_0$ семейство диффеоморфизмов окружности.

Если зависимость $f_s$ от $s$ гельдерова, то будем говорить 
о \emph{гельдеровом} косом произведении; наконец, если семейство $f_s$ 
является ограничением на $\Lambda_0$ гладко зависящего от $s\in D$
семейства, то такое семейство будем называть \emph{гладким}.

Гельдеровы косые произведения естественно возникают при ``выпрямлении''
ди\-на\-ми\-ки на частично-гиперболических инвариантных множествах.
Следующее предложение является связующим звеном между исследованием 
типичных диффеоморфизмов и теорией косых произведений. 
\begin{theorem}
\label{reduction}
Пусть задан гладкий диффеоморфизм $F$ из $D\times S^1$ в себя, являющийся 
косым произведением над отображением соленоида $T$, достаточно близкий 
к тождественному вдоль слоя. Тогда у любого диффеоморфизма $\widetilde{F}$, 
достаточно близкого к $F$, максимальный аттрактор $\widetilde{\Lambda}$ 
гомеоморфен $\Lambda_0\times S^1$, а ограничение 
$\widetilde{F}|_{\widetilde{\Lambda}}$ при этом го\-мео\-мор\-физ\-ме переходит
в \emph{гельдерово} косое произведение, близкое к $F|_{\Lambda_0\times S^1}$. 
\end{theorem}

Строгая формулировка и доказательство содержится в работах 
А.~С.~Городецкого~\cite{G05} и Hirsh, Pugh, Shub~\cite{HPS}.
Применение этой теоремы, позволяет перенести некоторые свой\-ства
косых произведений на диффеоморфизмы. Основной результат настоящей статьи
в терминах косых произведений формулируется следующим образом:

\begin{theorem}
\label{soft}
Существует открытая область $U$ в пространстве гёльдеровых косых
про\-из\-ве\-де\-ний над соленоидом Смейла-Вильямса $\Lambda_0$, содержащая гладкие 
косые произведения, сколь угодно близкие к тождественному по слою, 
такая, что любое ото\-бра\-же\-ние $f \in U$ имеет эргодическую меру
с полным носителем $\Lambda_0$ и нулевым показателем Ляпунова вдоль слоя.
\end{theorem}

Основной результат выводится из этих двух утверждений следующим образом:
рас\-сма\-три\-ва\-ет\-ся гладкое отображение компактного многообразия, обладающее
частично гиперболическим инвариантным множеством, гомеоморфным
прямому произведению соленоида на окружность $\Lambda_0 \times S^1$. 
Отображение строится таким образом, чтобы динамика
на инвариантном множестве являлась гладким косым произведением из области,
обес\-пе\-чен\-ной теоремой~\ref{soft}. Согласно теореме~\ref{reduction},
все диффеоморфизмы,  достаточно близкие к построенному, содержат
подмножество, динамика на котором сопряжена гёль\-де\-ро\-во\-му
косому произведению, близкому к исходному, а значит, принадлежащему
той же области. Тогда для этого косого произведения
найдётся инвариантная мера с нулевым показателем Ляпунова вдоль слоя.
Отсюда для соответствующей инвариантной меры возмущённого диффеоморфизма
один из показателей Ляпунова равен нулю.


\section{Гельдеровы косые произведения}
\label{SECTION_skew}

Рассмотрим соленоид Смейла-Вильямса, реализованный как отображение
полнотория. А именно, пусть $B$~--- единичный диск на плоскости
с комплексной координатой $u$, $S^1$~--- окружность, $D = S^1 \times B$.
Рассмотрим отображение $T: D \to D$:
$$
T: (\varphi, u) \to 
   (6\varphi\mod 1, \, \frac{1}{2} e^{2\pi i\varphi} + \frac{1}{100}u)
$$

Отображение $T$ имеет локально максимальный
ги\-пер\-бо\-ли\-чес\-кий ат\-трак\-тор~$\Lambda$, гомеоморфный $\SigmaSix$ --- 
пространству двусторонних пос\-ле\-до\-ва\-тель\-нос\-тей из символов $\{0, \dots, 5\}$,
в котором отождествлены точки $\{ \dots w_0 \dots w_k 5 5 5 \dots\}$ и
$\{ \dots w_0 \dots (w_k+1) 0 0 0 \dots \}$.
Ограничение $T|_{\Lambda}$ сопряжено сдвигу Бернулли
$\sigma: \SigmaSix \to \SigmaSix$, а отображениe~\eqref{solen_s1}, тем самым,
сопряжено отображению:
\begin{equation}
\label{skew}
G: \SigmaSixS \to \SigmaSixS, \quad (\omega, \varsloj) \to 
    (\sigma \omega, g_{\omega} (\varsloj)),
\end{equation}
где $g_{\omega}$~--- семейство диффеоморфизмов окружности, 
непрерывно в $\Diff^1$-норме зависящих от $\omega$
(в частности, гладкое вдоль слоёв).

Теорема~\ref{reduction} влечет выполнение для отображения $G$ 
{\bf условий $(L, C, \alpha)$-гельдеровости}:
$$
\forall \omega,\omega'\in \SigmaSix \qquad 
    d_{C^0}(g_\omega,g_{\omega'})< C \cdot d_{\widetilde{\Sigma}}
    (\omega,\omega')^{\alpha},
$$
$$
\forall \omega\in \SigmaSix \qquad
    \max_{x\in S^1}\max(g_{\omega}'(x),
    (g_{\omega}^{-1})'(x))<L.
$$
для некоторых констант $L > 1$, $C > 0$, $\alpha \in (0, 1)$.

\subsection{Управляемые косые произведения}

Введем обозначения:
$$
\begin{array}{rcl}
\bar g_m [\omega] & = &
    g_{\sigma^{m-1}\omega} \circ \dots \circ g_{\sigma\omega} \circ g_{\omega}, \\
\bar g_{-m} [\omega] & = &
    g^{-1}_{\sigma^{-m}\omega} \circ \dots \circ g^{-1}_{\sigma^{-1}\omega}, \\
\bar g_0 [\omega] & = & \id.
\end{array}
$$

Пусть $a$~--- произвольное конечное слово из символов $0, \dots, 5$.
Через $\{ \dots | a \dots \}$ мы обозначаем произвольную бесконечную
последовательность $\omega \in \SigmaSix$, у которой, на\-чи\-ная с
нулевого места, записано слово $a$. Аналогично вводим обозначения
$\{ \dots a | \dots \}$ и $\{ \dots a | b \dots \}$.
Латинская $w$ и греческая $\omega$ обозначают слова конечной
и бесконечной длины соответственно.

\begin{definition}
\label{upr_prop}
Пусть косое произведение~$G$ является $(L,C,\alpha)$-гельдеровым.
Скажем, что оно обладает свойством\par
\begin{itemize}
\item
\emph{растяжения,} если существуют $\nu > 1$, $\delta_1>0$, такие,
что для произвольного интервала $I \subset S^1$, $|I| < \delta_1$:
\begin{equation}
\label{prop_j1}
\exists j_1 \neq 5:  \, \forall \omega
    = \{ \dots | j_1 \beta \dots \}, \beta \neq 5, \quad \forall x\in I
    \quad (Dg_{\omega})(x) > \nu
\end{equation}
\par
\item
\emph{обратного растяжения,} если существуют $\nu > 1$,
$\delta_1>0$, такие, что для про\-из\-воль\-но\-го интервала $I \subset
S^1$, $|I| < \delta_1$:
\begin{equation}
\label{prop_j2}
\exists j_2 \neq 5:  \, \forall \omega
    = \{ \dots j_2 | \beta \dots \}, \beta \neq 5, \quad \forall x\in I
    \quad (Dg^{-1}_{\sigma^{-1} \omega})(x) > \nu;
\end{equation}
\par
\item
\emph{наличия $\delta_2$-поворота,} если для любой
последовательности $\omega = \{ \dots | 0 \beta \dots \}$, $\beta \neq 5$,
верно:
\begin{equation}
\label{povorot} d_{C^0} ( g_{\omega}, H_{\delta_2} ) <
\frac{\delta_2^2}{40};
\end{equation}
здесь через $H_{\delta_2}$ обозначен поворот на угол~$\delta_2$.
\par
\item
\emph{наличия слабопритягивающей орбиты,} если существует
притягивающая периодическая орбита $X$, для показателя Ляпунова
вдоль слоя которой, обозначаемого $\lambda(X)<0$, выполнено
\begin{equation}\label{X_1}
\lambda(X) + \ln\nu > 0
\end{equation}
\par
\item
\emph{$\gamma$-предсказуемости траекторий,} $\gamma>0$, если для
любой точки $\varsloj \in S^1$, для любого натурального $m$, для
любого конечного слова $w^{*} = w_{-m} \dots w_{-1} | w_0 \dots
w_{m-1}$ вы\-пол\-ня\-ет\-ся
\begin{eqnarray}
\label{prop_gamma} \diam \{ \bar g_{\pm m}[\omega] (\varsloj) | \omega =
\{ \dots w^{*} \dots \} \}
    & < & \gamma
\end{eqnarray}
\end{itemize}
\end{definition}

\begin{definition}
$(L,C,\alpha)$-гельдерово косое произведение называется \emph{управляемым},
если оно обладает всеми свойствами определения~\ref{upr_prop},
причём константы могут быть выбраны удовлетворяющими
\emph{условию согласованности констант}. Последнее состоит в том, что
\begin{equation}
\label{gammasize} \gamma < \frac{\delta_2}{40}, \quad
\delta_1>3\delta_2
\end{equation}
и
$$
\alpha > \log_2 L.
$$
\end{definition}

\begin{lemma}
\label{construct}
Для любых заданных $L>1$, $C>0$ и $\alpha\in(0,1)$
множество управляемых систем непусто, открыто и содержит отображения,
сколь угодно близкие к тождественному вдоль слоя.
\end{lemma} 

\begin{pf}
Действительно, в произвольной окрестности $W \subset \Diff^1(S^1)$  
тождественного отображения окружности можно выбрать три гиперболических
диффеоморфизма $g_{1,2,3}$ с одним аттрактором и одним репеллером,
обеспечивающие свойства рас\-тя\-же\-ния и обратного растяжения,
поворот на малый угол $g_0$, диффеоморфизм Морса-Смейла
$g_4$ с притягивающей периодической орбитой и тождественное
отображение $g_5 := id$.
Определим отображения семейства $g_{\omega} := g_i$ на словах вида
$\{ \dots | i \beta \dots \}$, где $\beta \neq 5$, и продолжим семейство
гладким образом на оставшиеся слова базы. Отображения продолжения
можно выбрать принадлежащими окрестности $W$. Непустота доказана.

В отношении всех свойств определения~\ref{upr_prop}, кроме
свойства предсказуемости тра\-ек\-то\-рий, ясно, что они $C^1$-устойчивы.
Следующая лемма, принадлежащая А.~С.~Городецкому, выводит предсказуемость
траекторий из условия на показатель Гёльдера и скорость рас\-тя\-же\-ния.

\begin{lemma}[{\cite[Lemma~3.1]{GI01}}]\label{sloj}
Пусть заданы $L$, $C$ и $\alpha$, такие, что выполнено
$\alpha > \log_2 L$~--- второе из условий согласованности констант.
Тогда существует такое $K = K (L, C, \alpha)$, что для любой
$(L,C,\alpha)$-системы, $\delta$-близкой к $S$-ступенчатой,
\begin{equation}
d_{C_0} (\bar g_{\pm m} [\omega], \bar g_{\pm m} [\omega']) \leq
\gamma := K\delta^{\beta},
\end{equation}
где $\beta = 1 - \frac{\ln L}{\ln 2^{\alpha}}$.
\end{lemma}

Таким образом, неравенство $K\delta^{\beta}< \frac{\delta_2}{40}$ на
расстояние $\delta$ от управляемой системы до её возмущения влечет
выполнение условия предсказуемости траекторий для последнего.
\end{pf}

\subsection{Построение периодических орбит}

\begin{definition}
\label{goodpoint}
Пусть задана периодическая орбита $X$ отображения $G$,
период $X$ равен $P$, и \mbox{$\varepsilon > 0$}~-- фиксировано.
Точка $y$ называется \emph{$(\varepsilon, P)$-близкой}
к орбите $X$, если найдется $x \in X$ такое, что
$$
\forall l = 0, 1 \dots P-1 \qquad d(G^l(x), G^l(y)) < \varepsilon.
$$
\end{definition}

\begin{lemma}[Основная лемма]
\label{mainlemma}
Пусть косое произведение $G$ вида~(\ref{skew})
обладает свойством управляемости и $D = D(G) > 0$~--- некоторая константа. 
Пусть $X$~--- про\-из\-воль\-ная периодическая орбита $G$ периода $P$ 
с мультипликатором по слою $0 < \theta < 1$, для показателя Ляпунова
вдоль слоя $\lambda := \frac{\ln \theta} {P}$ которой выполняется:
$$
\lambda + \ln \nu > D.
$$
Пусть $U$~-- произвольная окрестность в $\SigmaSixS$
Тогда для любого $\varepsilon > 0$ существует периодическая орбита
$Y$ косого произведения $G$ с периодом $P' > 2P$ и показателем
Ляпунова вдоль слоя $\lambda' < 0$ такая, что:

\begin{enumerate}

\item орбита $Y$ пересекает окрестность $U$;

\item $|\lambda'| < C |\lambda| $, где $C = C(G)$~--- глобальная константа,
зависящая только от косого произведения $G$, $0 < C < 1$;

\item $\lambda' + \ln\nu > D$;

\item\label{i:projection} существует ${\widetilde Y} \subset Y$
и проекция $\pi: {\widetilde Y} \to X$ такие, что:

\begin{enumerate}

\item \label{2a}
все точки множества~${\widetilde Y}$ являются
$(\varepsilon,P)$-близкими к орбите $X$, причем в
определении~\ref{goodpoint} можно взять $x = \pi (y)$;

\item \label{2b}
доля $\varkappa := \frac { \# {\widetilde Y} }{ \# Y }$ точек, в
которых определена проекция $\pi$, оценивается как:
$$
\varkappa \geq 1 - \frac{3|\lambda|}{\ln L};
$$

\item \label{2c}
количество элементов прообраза $\pi^{-1} (x)$ одинаково для всех $x
\in X$.

\end{enumerate}

\end{enumerate}

\end{lemma}

\begin{pf}
Периодическая орбита косого произведения задается своей начальной
точкой $(\omega, \varsloj)$, $\omega \in \SigmaSix$, $\varsloj \in
S^1$, $\omega = (w)$~--- периодическая последовательность, $w = (w_0
\ldots w_{P-1})$~--- ее период и верно:
$$
\sigma^P \omega = \omega \qquad \bar g_P [\omega] (\varsloj) =
\varsloj.
$$

Используя свойства управляемости и зафиксировав интервал $J$
на центральном слое, мы подбираем серию слов $w'(k)$ в базе,
длина которых возрастает с ростом натурального параметра $k$.
Соответствующие отображения
косого произведения будут сжимать $J$ в себя, гарантируя наличие
притягивающей вдоль слоя неподвижной точки, а их про\-из\-вод\-ная на $J$
будет равномерно ограничена по $k$. Увеличивая $k$, можно добиться
сколь угодно близкого к нулю отрицательного показателя Ляпунова.

Попутно, строя новую орбиту $Y$, мы обеспечиваем ее прохождение через 
окрестность $U$ и близость точек новой и старой орбит 
(свойства {\it 4a}--{\it 4c}). Это позволит нам получить свойства предельной
меры: весь аттрактор в качестве носителя и близость
временных и пространственных средних.

Техника доказательства подробно описана в~\cite{KN06}.
\end{pf}

\subsection{Эргодичность и нулевые показатели Ляпунова}
\label{SECTION_proof}

\begin{lemma}\label{l:zero}
Для любой управляемой системы найдётся эргодическая инвариантная
мера с полным носителем и нулевым показателем Ляпунова вдоль слоя.
\end{lemma}
\begin{pf}
Выберем счётное семейство окрестностей $U_i$ таких, что каждая точка
пространства покрыта окрестностью сколь угодно малого диаметра.

Воспользовавшись леммой~\ref{mainlemma}, мы можем построить
последовательность при\-тя\-ги\-ва\-ю\-щих вдоль слоя периодических орбит $X_i$,
начинающуюся со слабопритягивающей ор\-би\-ты (см.~\eqref{X_1}),
каждая орбита $X_i$ пересекается с окрестностью $U_i$.
Показатели Ляпунова для этих орбит экспоненциально стремятся к нулю,
причём каждая следующая орбита б\'ольшую часть времени проводит около
предыдущей.

Рассмотрим последовательность атомарных мер, равномерно
распределённых на этих орбитах. Из условия на ``похожесть'' орбит
(см.~заключение~\ref{i:projection} леммы~\ref{mainlemma}) с помощью
эр\-го\-ди\-чес\-кой теоремы Биркгофа-Хинчина выводится, что любая
предельная точка пос\-тро\-ен\-ной последовательности будет эргодической
инвариантной мерой; также несложно проверяется неатомарность
предельной меры. Поскольку пространство мер на $\SigmaSixS$
$*$-слабо компактно, мы можем выделить из этой последовательности
сходящуюся~--- в силу вышесказанного к эргодической инвариантной
мере~--- подпоследовательность. Эта предельная мера и будет искомой.

Действительно, показатель Ляпунова для эргодической инвариантной
меры выражается как интеграл по этой мере от непрерывной функции~---
производной отображения вдоль слоя. Поэтому показатель Ляпунова для
предельной меры равен пределу показателей Ляпунова, то есть нулю.

С другой стороны, если орбита $X_i$ перескает $U_i$, то пересечение
есть и у всех орбит $X_j$, $j > i$, причем доля общих точек 
не стремится к нулю. Тем самым, предельная мера $U_i$ не равна нулю,
что обеспечивает совпадение носителя меры со всем пространством.

Эти рассуждения строго проведены в работе~\cite{GIKN} (см. леммы~1
и~2) для частного случая ступенчатых косых произведений
над подковой Смейла и дословно переносятся на общий случай.
\end{pf}

\subsection{Доказательство теоремы~\ref{soft}}

Построенная при доказательстве леммы~\ref{construct} система
является не только гельдеровой, но и гладкой, а все близкие --- управляемыми
гельдеровыми. В силу лемм~\ref{mainlemma} и~\ref{l:zero} для них
найдется эргодическая инвариантная мера с нулевым показателем Ляпунова
вдоль слоя. Теорема доказана.

\end{document}